\documentclass[a4paper,11pt]{article}

\usepackage{amssymb,amsmath,amsthm}
\usepackage{mathrsfs}
\usepackage{authblk}
\usepackage{indentfirst}
\usepackage[numbers,sort&compress]{natbib}
\usepackage[english]{babel}
\usepackage{caption}
\usepackage{graphicx}
\usepackage{hyperref}
\usepackage{bookmark}
\usepackage{tikz}
\usepackage{helvet}
\usepackage[EULERGREEK]{sansmath}
\usepackage{enumitem}
\usepackage[section]{placeins}
\usetikzlibrary{arrows.meta,calc,positioning}

\newtheorem{theorem}{Theorem}[section]
\newtheorem{corollary}[theorem]{Corollary}
\newtheorem{lemma}[theorem]{Lemma}
\newtheorem{proposition}[theorem]{Proposition}
\theoremstyle{definition}
\newtheorem{definition}[theorem]{Definition}
\newtheorem{example}[theorem]{Example}
\newtheorem{remark}[theorem]{Remark}

\DeclareMathOperator{\im}{im}

\DeclareMathOperator{\Spec}{Spec}

\newcommand{\R}{\mathbb R}
\newcommand{\A}{\mathcal A}
\newcommand{\LH}{\mathrm{LH}}
\newcommand{\loc}{\mathrm{loc}}
\newcommand{\id}{\mathrm{id}}
\graphicspath{{figures/}}
\makeatletter
\renewcommand*\l@subsection{\@dottedtocline{2}{1.5em}{2.8em}}
\makeatother
\hypersetup{
  colorlinks=false,
  pdfborder={0 0 1},
  linkbordercolor={1 0 0},
  citebordercolor={0 1 0},
  urlbordercolor={0 1 1},
  pdftitle={Local path Laplacian},
  pdfauthor={Hongsong Feng and Jian Liu}
}

\title{Local path Laplacian for directed high-order interactions}

\author[1]{Jian Liu}
\author[2]{Hongsong Feng\thanks{Corresponding author: hfeng2@charlotte.edu}}
\affil[1]{Mathematical Science Research Center, Chongqing University of Technology, Chongqing 400054, China}
\affil[2]{Department of Mathematics and Statistics, University of North Carolina at Charlotte, Charlotte, NC 28223, USA}

\date{}

\begin{document}
\maketitle

\paragraph{Abstract}

Path homology, introduced by Grigor'yan, Lin, Muranov and Yau (GLMY) to characterize the higher-order topology of directed networks, has found applications across a wide range of fields.
Building on this theory, we introduce a local path Laplacian for spectral analysis near a marked set of vertices in a digraph. This self-adjoint, positive semidefinite operator captures the local structure of directed paths. Its kernel is naturally isomorphic to local path homology, and its spectrum is an intrinsic geometric invariant of the digraph with its marked set.
We prove that a finite neighborhood of the marked set determines the local path Laplacian up to unitary equivalence. This result has two practical consequences. First, computing the local spectrum requires only information from this neighborhood and can therefore be substantially faster than computing the global path spectrum when the neighborhood is small relative to the network. Second, local spectra at different marked sets can be computed independently and in parallel. When the network is extended, only marked sets whose defining neighborhoods change require recomputation.
We evaluate how well local path spectra distinguish isomorphism classes of five-vertex digraphs without reciprocal edges. For neuron matching between two adult Caenorhabditis elegans connectomes, adding the local spectrum in dimension one to the in- and out-degrees raised the top-five hit rate from 34.06 percent to 42.22 percent.

\paragraph{Keywords}
Digraphs; path homology; quotient inner products; local path Laplacian; local spectra.
\footnotetext[1]{{\bf 2020 Mathematics Subject Classification.} Primary 55N35; Secondary 05C20, 05C50.}

\begingroup
\linespread{1}\small
\tableofcontents
\endgroup
\clearpage

\section{Introduction}

Graphs and digraphs are among the most widely used models of complex networks~\cite{BattistonHigherOrder}. Classical graph invariants include the degree sequence, the chromatic polynomial, and the spectrum of the graph Laplacian~\cite{GodsilRoyle}. These invariants describe the network through its pairwise connections. For example, three separate pairwise interactions and one interaction involving all three vertices can give the same triangle. A classical way to model higher-order structure is to form the flag complex, also called the clique complex, by taking each clique of an undirected graph as a simplex. For a digraph, the directed flag complex uses ordered cliques in which each vertex has an arrow to every later vertex. The homology of clique complexes has been used to detect cavities in neuronal networks~\cite{ReimannCliques}, and combinatorial Laplacians on clique complexes have been used to identify cyclic inconsistencies in ranking data~\cite{JiangHodgeRanking}.

Grigor'yan, Lin, Muranov and Yau (GLMY) introduced path complexes and path homology to study the higher-order structure of digraphs~\cite{GLMY}. Their construction associates to a digraph $G$ the path complex $P(G)$ of directed paths. Its homology is called \emph{path homology}. Every simplicial complex with ordered vertices gives a path complex by reading each simplex in increasing order~\cite{GLMY}. In topological data analysis, persistent path homology and persistent path Laplacians track features of directed networks across scales~\cite{ChowdhuryMemoliPath,WangWeiPathLaplacian}. Applications of path topology include the analysis of molecular and materials data~\cite{ChenLiuPathTopology} and the design of network architectures for reservoir computing~\cite{ChenWangLin2026}.

The usual path homology and path Laplacian describe a directed network as a whole. To study the structure near a prescribed vertex or set of vertices, one needs an invariant that also depends on the choice of those vertices. There is a computational difficulty as well. The number of allowed paths can grow rapidly with dimension, and computing the invariant chain spaces requires solving linear boundary constraints~\cite{ChowdhuryMemoliPath,DeyLiWangPath,WangWeiPathLaplacian}. A further difficulty arises when the network evolves. For example, under the standard path inner product, adding a single arrow to a loopless digraph increases the trace of its path Laplacian in dimension $0$ by two. A model that uses global spectral features may therefore need to update its inputs after a local change. Repeated updates of graph representations are a practical concern in learning on dynamic networks~\cite{KazemiDynamicGraphs}.

Motivated by these questions and by the local simplicial Laplacian~\cite{LiuFengLiu}, we introduce the \emph{local path Laplacian} at a marked vertex set of a regular path complex. The local path Laplacian is a self-adjoint and positive semidefinite operator, and its kernel is naturally isomorphic to the local path homology (Theorems~\ref{thm:laplacian-properties} and~\ref{thm:hodge}). Theorem~\ref{thm:spectral-locality} shows that this operator is unitarily equivalent to the local path Laplacian on an induced neighborhood path complex. When this neighborhood is small, the matrix construction can use far fewer paths than a computation on the full complex. An extension of the network that preserves the induced neighborhood preserves every local eigenvalue and its multiplicity. For a new vertex whose neighborhood lies in the added subgraph, that subgraph alone determines its local spectrum.

Furthermore, we give explicit matrix computations for some basic examples. We also show that the computational complexity of the local path Laplacian is substantially lower than that of the global path Laplacian for large networks. As an experiment, we compute local spectra in dimensions $0$ and $1$ for all $582$ isomorphism types of five-vertex digraphs without parallel edges. Among the $394$ pairs of nonisomorphic graphs with the same multiset of in-degree and out-degree pairs, recording the local spectrum in dimension $1$ alongside the degree pair at each vertex distinguishes $358$ pairs. In two adult \emph{Caenorhabditis elegans} connectomes, adding the local spectrum raises the top-five hit rate for matching neurons across animals from $34.06\%$ to $42.22\%$ (Section~\ref{sec:adult-neuron-matching}).

This work is organized as follows.
Section~\ref{sec:definitions} defines the local path Laplacian and gives its matrix representation.
Section~\ref{sec:properties} establishes the Hodge decomposition and shows that a finite neighborhood determines the local spectrum.
Section~\ref{sec:examples} illustrates how branching and vertex marking affect the spectrum.
Section~\ref{sec:experiments} examines graph discrimination, neuron matching across connectomes, and the cost of computing local and global spectra.

\section{Path complexes and the local path Laplacian}\label{sec:definitions}

In this section, we recall the basic definitions of path complexes and path homology. We then equip the quotient chain complex with an inner product structure and define the local path Laplacian. Finally, we derive a matrix representation of the local path Laplacian and use it to compute its eigenvalues.

\subsection{Regular path chains}

All vector spaces are over $\R$. Let $V$ be a finite nonempty vertex set.

\begin{definition}\label{def:paths}
An \textit{elementary $p$-path} on $V$ is a sequence $x_0\cdots x_p$, denoted by $e_{x_0\cdots x_p}$. It is \textit{regular} if $x_i\ne x_{i+1}$ for every $0\le i<p$. Otherwise, it is \textit{nonregular}.
\end{definition}

\begin{definition}\label{def:path-complex}
A \textit{path complex} $P$ on $V$ is a collection of paths with the following property. If $x_0\cdots x_p\in P$ and $p\ge1$, then $x_0\cdots x_{p-1}\in P$ and $x_1\cdots x_p\in P$. The paths in $P$ are called \textit{allowed}.
If all paths are regular, we call $P$ a \textit{regular path complex}.
\end{definition}

Let $\mathcal R_p(V)$ be the vector space spanned by the regular elementary $p$-paths on $V$. The regular boundary on $\mathcal R_*(V)$ is given by
\[
 \partial_p e_{x_0\cdots x_p} =\sum_{i=0}^p(-1)^i e_{x_0\cdots\widehat{x_i}\cdots x_p}
\]
for $p\ge1$, with nonregular faces interpreted as zero. Set $\partial_0=0$. These maps satisfy $\partial_{p-1}\partial_p=0$ for $p\ge2$.

For a regular path complex $P$, let $\A_p(P)$ be the vector space generated by the $p$-paths in $P$. Then $\A_p(P)$ is a subspace of $\mathcal R_p(V)$. The \textit{$\partial$-invariant chain complex} is
\[
 \Omega_p(P)=\{c\in\A_p(P):\partial_pc\in\A_{p-1}(P)\},
 \quad p\ge1,
\]
and $\Omega_0(P)=\A_0(P)$. The \textit{path homology} of $P$ is the homology of the chain complex $\Omega_\ast(P)$. Path homology is a topological invariant of path complexes that captures higher-order features of directed systems. Its applications include the study of directed networks, molecular and materials systems, and reservoir computing \cite{ChowdhuryMemoliPath,ChenLiuPathTopology,ChenWangLin2026}.

Throughout this paper, digraphs are finite and have no loops. The path complex $P(G)$ of a digraph $G=(V,E)$ consists of all vertex sequences $x_0\cdots x_p$ such that $(x_{i-1},x_i)\in E$ for $1\le i\le p$. A path may revisit a vertex. The path homology of $G$ is the path homology of $P(G)$.

\subsection{Path inner products and quotient representatives}

For a subset $U\subseteq V$, let $P|_U$ be the induced path complex of allowed paths whose vertices lie in $U$. Given a marked set $T\subseteq V$, write $P-T=P|_{V\setminus T}$. In the digraph case, this is the path complex of the induced digraph $G-T$.

\begin{definition}\label{def:local-complex}
Let $\Omega_*(P)$ and $\Omega_*(P-T)$ denote the $\partial$-invariant chain complexes of $P$ and $P-T$, respectively. The \textit{local path complex} at $T$ is their quotient, with chain groups and differential given by
\[
 Q_p(P,T)=\frac{\Omega_p(P)}{\Omega_p(P-T)}, \quad d_p[c]=[\partial_pc].
\]
The homology of $Q_*(P,T)$ is called the \textit{local homology} at $T$ and is denoted by $\LH_*(P,T)$ \cite{GJMY,LocalHomology}.
\end{definition}

For a single marked vertex $v$, we write $Q_*(P,v)$ and $\LH_*(P,v)$. For a digraph, we use $Q_*(G,T)$ and $\LH_*(G,T)$. In particular, $Q_*(P,\varnothing)=0$ and $Q_*(P,V)=\Omega_*(P)$. The \textit{local Betti number} in dimension $p$ is $\beta_p(P,T)=\dim\LH_p(P,T)$.

The invariant chains of the deletion complex satisfy
\[
 \Omega_p(P-T)=\Omega_p(P)\cap\mathcal R_p(V\setminus T).
\]

We take the elementary allowed $p$-paths of $P$ as an orthonormal basis of $\A_p(P)$. The subspace $\Omega_p(P)$ inherits this inner product. Let
\[
 W_p(P,T)=\Omega_p(P)\cap \Omega_p(P-T)^\perp.
\]
The subspace $W_p(P,T)$ inherits the inner product on $\Omega_p(P)$ and is naturally isomorphic to the quotient space $Q_p(P,T)$.

\begin{definition}\label{def:metric}
The \textit{quotient inner product} on $Q_p(P,T)$ is the unique inner product that makes the map below an isometry.
\[
 W_p(P,T)\longrightarrow Q_p(P,T),\quad w\longmapsto[w].
\]
\end{definition}

Let $\Pi_p:\Omega_p(P)\to W_p(P,T)$ be the orthogonal projection. For $p\ge1$, define $d_p:W_p(P,T)\to W_{p-1}(P,T)$ by
\[
 d_p=\Pi_{p-1}\partial_p\big|_{W_p(P,T)},
\]
and set $W_{-1}(P,T)=0$ and $d_0=0$. 
\begin{proposition}\label{prop:projected-boundary}
The pair $(W_*(P,T),d_*)$ is a chain complex. The maps
\[
 q_p:W_p(P,T)\longrightarrow Q_p(P,T),\quad q_p(w)=[w],
\]
give a chain isomorphism
\[
 q_*:(W_*(P,T),d_*)\xrightarrow{\cong}
 \left(Q_\ast(P,T),\overline\partial_*\right),
\]
where $\overline\partial_p[c]=[\partial_pc]$ is the quotient differential.
\end{proposition}
\begin{proof}
By the preceding construction, each $q_p$ is an isometric linear isomorphism with inverse $[c]\mapsto\Pi_pc$. Since $\Omega_*(P-T)$ is a subcomplex, $\overline\partial_*$ is well defined and satisfies $\overline\partial_{p-1}\overline\partial_p=0$. For any $w\in W_p(P,T)$ and $p\ge1$, we have
\[
 \partial_pw-\Pi_{p-1}\partial_pw\in\Omega_{p-1}(P-T).
\]
It follows that
\[
 q_{p-1}d_pw=[\Pi_{p-1}\partial_pw]
 =[\partial_pw]=\overline\partial_pq_pw.
\]
Thus $W_*(P,T)$ is a chain complex and $q_*$ is the asserted chain isomorphism.
\end{proof}

Under this isomorphism, each class $[c]\in Q_p(P,T)$ has the unique representative $\Pi_pc$ of minimum norm. Thus, we have
\[
 \langle[c],[c']\rangle_{Q_p}
 =\langle\Pi_pc,\Pi_pc'\rangle,
 \quad
 \|[c]\|_{Q_p}^2
 =\min_{k\in\Omega_p(P-T)}\|c+k\|^2.
\]
Indeed, every representative of $[c]$ has the form
$\Pi_pc+k$ with $k\in\Omega_p(P-T)$, and
\[
 \|\Pi_pc+k\|^2=\|\Pi_pc\|^2+\|k\|^2.
\]
In particular, the norm of an arbitrary lift of $[c]$ need not equal its quotient norm.

\subsection{The local path Laplacian}

\begin{definition}\label{def:laplacian}
Let $d_p^*:Q_{p-1}(P,T)\to Q_p(P,T)$ denote the adjoint of $d_p$ with respect to the quotient inner products. The \textit{local path Laplacian} in dimension $p$ is
\[
 \Delta_p^{\loc}(P,T)=d_p^*d_p+d_{p+1}d_{p+1}^*.
\]
The summand $d_p^*d_p$ is the down Laplacian, while the summand $d_{p+1}d_{p+1}^*$ is the up Laplacian. 
\end{definition}

When $T=V$, the local path Laplacian reduces to the usual path Laplacian on $\Omega_*(P)$. When $P$ comes from an ordered simplicial complex and $T$ consists of a single vertex, this construction gives the local simplicial Laplacian at that vertex \cite{LiuFengLiu}. For a digraph $G$, we write $\Delta_p^{\loc}(G,T)$ for $\Delta_p^{\loc}(P(G),T)$.

\begin{remark}\label{rem:compression}
The local operator is generally not obtained by compressing the absolute path Laplacian. Let $\Delta_p^\Omega$ be the Laplacian on $\Omega_*(P)$ and set
\[
 t_p=(\id-\Pi_{p-1})\partial_p\big|_{W_p}:W_p\longrightarrow \Omega_{p-1}(P-T).
\]
Then the decomposition $\Omega_p(P)=\Omega_p(P-T)\oplus W_p(P,T)$ gives
\[
 \Pi_p\Delta_p^\Omega\big|_{W_p}
 =\Delta_p^{\loc}+t_p^*t_p.
\]
Indeed, the boundary has the block form
$\left(\begin{smallmatrix}a_p&t_p\\0&d_p\end{smallmatrix}\right)$
because $\Omega_*(P-T)$ is a subcomplex. The discarded boundary component $t_p$ accounts for the correction term $t_p^*t_p$.
\end{remark}

\subsection{Matrix representations and computation}\label{sec:matrices}

Choose an orthonormal basis of each $Q_p(P,T)$, and let $B_p$ be the matrix of $d_p$ in these bases. Then the local path Laplacian has matrix
\[
 L_p^{\loc}=B_p^{\mathsf T}B_p+B_{p+1}B_{p+1}^{\mathsf T}.
\]
In particular, $L_0^{\loc}=B_1B_1^{\mathsf T}$.

More generally, choose a basis in each dimension. Let $M_p$ be the Gram matrix of the basis of $Q_p$, and let $D_p$ be the matrix of $d_p$. The matrices of the adjoint and the local Laplacian are then
\[
 \begin{aligned}
 D_p^{\dagger}&=M_p^{-1}D_p^{\mathsf T}M_{p-1},\\
 L_p^{\loc}
 &=M_p^{-1}D_p^{\mathsf T}M_{p-1}D_p
   +D_{p+1}M_{p+1}^{-1}D_{p+1}^{\mathsf T}M_p.
 \end{aligned}
\]
When the basis is orthonormal, the Gram matrix is the identity matrix.

Fix an ordered basis of elementary allowed paths in each dimension, identifying $\A_q(P)$ with $\R^{n_q}$. For $q\ge1$, write the boundary matrix as
\[
 \partial_q\big|_{\A_q(P)}=
 \begin{pmatrix}E_q\\F_q\end{pmatrix}.
\]
The rows of $E_q$ and $F_q$ correspond to allowed and non-allowed regular $(q-1)$-paths, respectively. A chain with coordinate vector $x$ is $\partial$-invariant precisely when $F_qx=0$. Thus $\Omega_q(P)=\ker F_q$.

Let the columns of $C_q$ form a basis of $\Omega_q(P-T)$ in allowed-path coordinates, with zero entries on paths meeting $T$. Then the orthogonal complement of $\Omega_q(P-T)$ is the solution space of $C_q^{\mathsf T} x = 0$. Thus, we have
\[
 W_q(P,T)=\ker F_q\cap\ker C_q^{\mathsf T}
 =\ker\begin{pmatrix}F_q\\C_q^{\mathsf T}\end{pmatrix}.
\]
Choose $U_q\in\R^{n_q\times m_q}$ with orthonormal columns spanning $ W_q(P,T)$, where $m_q=\dim W_q(P,T)$. Then its columns also represent an orthonormal basis of $Q_q(P,T)$.

For $q\ge1$, the identity $F_qU_q=0$ gives $\partial_qU_q=E_qU_q$ in allowed-path coordinates. Note that the orthogonal projection $\Pi_{q-1}$ is represented by $U_{q-1}U_{q-1}^{\mathsf T}$. If $B_q$ denotes the matrix of $d_q$ in the chosen bases, then we have
\[
 U_{q-1}B_q
 =d_qU_q
 =\Pi_{q-1}\partial_qU_q
 =U_{q-1}U_{q-1}^{\mathsf T}E_qU_q.
\]
Multiplying on the left by $U_{q-1}^{\mathsf T}$ and using $U_{q-1}^{\mathsf T}U_{q-1}=I_{m_{q-1}}$, we obtain
\[
 B_q=U_{q-1}^{\mathsf T}E_qU_q\in\R^{m_{q-1}\times m_q}.
\]
The matrices $B_p$ and $B_{p+1}$ determine $L_p^{\loc}$. Theorem~\ref{thm:spectral-locality} shows that the computation requires only vertices in a bounded neighborhood of $T$.

%
%

\begin{figure}[htbp]
\centering
\includegraphics[width=\textwidth]{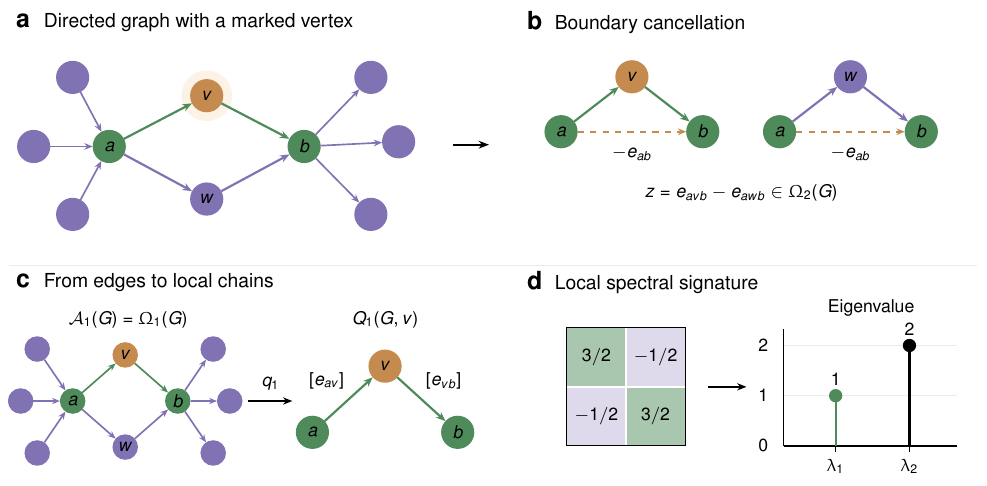}
\caption{\textbf{From a directed graph to local spectral features.}  Green arrows form an orthonormal basis of the local $1$-chain space, while purple arrows avoid the marked set.
\textbf{(a)} A directed graph with marked vertex $v$ in orange. Green and purple edges form the two paths $a\to v\to b$ and $a\to w\to b$.
\textbf{(b)} Each path has the non-allowed boundary term $-e_{ab}$. These terms cancel in $z=e_{avb}-e_{awb}$, giving a $\partial$-invariant $2$-chain.
\textbf{(c)} In dimension $1$, every allowed edge is $\partial$-invariant. The edge-basis diagrams show the quotient map $q_1$: edges avoiding $v$ map to zero, while $[e_{av}]$ and $[e_{vb}]$ form an orthonormal basis of $Q_1(G,v)$.
\textbf{(d)} The local Laplacian matrix and its eigenvalues $1$ and $2$.}
\label{fig:network-workflow}
\end{figure}

Figure~\ref{fig:network-workflow} illustrates the construction on a directed graph with a marked vertex $v$. The two paths from $a$ to $b$ have a common non-allowed boundary term $-e_{ab}$. Their difference is a $\partial$-invariant chain whose class contributes to the local up Laplacian.

\section{Related properties}\label{sec:properties}

From now on, we assume that $P$ is a path complex on $V$ and $T$ is a marked subset of $V$.

\subsection{Hodge decomposition and local homology}

\begin{theorem}\label{thm:laplacian-properties}
For any $p\ge0$, the operator $\Delta_p^{\loc}(P,T)$ is self-adjoint and positive semidefinite. Moreover, we have
\[
 \langle\Delta_p^{\loc}x,x\rangle
 =\|d_px\|^2+\|d_{p+1}^*x\|^2,
 \quad
 \ker\Delta_p^{\loc}=\ker d_p\cap\ker d_{p+1}^*.
\]
\end{theorem}
\begin{proof}
Both summands of $\Delta_p^{\loc}$ are self-adjoint. By the definition of the adjoint,
\[
 \langle\Delta_p^{\loc}x,x\rangle
 =\langle d_p^*d_px,x\rangle
  +\langle d_{p+1}d_{p+1}^*x,x\rangle
 =\|d_px\|^2+\|d_{p+1}^*x\|^2.
\]
The right-hand side is nonnegative, so $\Delta_p^{\loc}$ is positive semidefinite. If $\Delta_p^{\loc}x=0$, then $d_px=0$ and $d_{p+1}^*x=0$. Conversely, these two equalities imply $\Delta_p^{\loc}x=0$.
\end{proof}

\begin{theorem}\label{thm:hodge}
For any $p\ge0$, there is an orthogonal decomposition of vector spaces
\[
 Q_p(P,T)
 =\im d_{p+1}\mathbin{\oplus}
 \ker\Delta_p^{\loc}\mathbin{\oplus}
 \im d_p^*.
\]
Moreover, we have an isomorphism
\[
 \ker\Delta_p^{\loc}(P,T)\xrightarrow{\cong}\LH_p(P,T).
\]
In particular, the multiplicity of the zero eigenvalue equals the local Betti number $\beta_p(P,T)$.
\end{theorem}
\begin{proof}
The spaces $\im d_{p+1}$ and $\im d_p^*$ are orthogonal because $d_pd_{p+1}=0$. Their common orthogonal complement is $\ker d_{p+1}^*\cap\ker d_p$, which equals $\ker\Delta_p^{\loc}$ by Theorem~\ref{thm:laplacian-properties}. This proves the displayed decomposition.

A cycle is orthogonal to $\im d_p^*$ and therefore has a unique expression $z=d_{p+1}y+h$, where $h$ is harmonic. Hence, each homology class has exactly one harmonic representative. The claimed isomorphism follows.
\end{proof}

One may assign a positive squared norm $w_p(\alpha)$ to each elementary allowed path $\alpha$. In this case, Theorem~\ref{thm:hodge} also holds and the nullity does not depend on the positive weights. In this work, we always set $w_p(\alpha)=1$.

\subsection{Unitary equivalence}

\begin{proposition}\label{prop:isometry}
Let $F_*:Q_*(P,T)\to Q_*(P',T')$ be an isometric chain isomorphism. Then $\Delta_p^{\loc}(P,T)$ and $\Delta_p^{\loc}(P',T')$ are unitarily equivalent. More precisely, for any $p\ge0$, one has
\[
 \Delta_p^{\loc}(P',T')
 =F_p\Delta_p^{\loc}(P,T)F_p^*.
\]
\end{proposition}
\begin{proof}
Since $F_p$ is an isometric isomorphism, $F_p^{-1}=F_p^*$. Let $d_p$ and $d'_p$ denote the boundary operators on $Q_*(P,T)$ and $Q_*(P',T')$, respectively. For $p\ge1$, the chain map identity gives $d_p'=F_{p-1}d_pF_p^*$. Taking adjoints gives
\[
 \begin{aligned}
 \Delta_p^{\loc}(P',T')
 &=d_p'^*d_p'+d_{p+1}'d_{p+1}'^*\\
 &=F_p\bigl(d_p^*d_p+d_{p+1}d_{p+1}^*\bigr)F_p^*\\
 &=F_p\Delta_p^{\loc}(P,T)F_p^*.
 \end{aligned}
\]
For $p=0$, the same calculation applies with $d_0=d'_0=0$.
\end{proof}

\begin{corollary}\label{cor:reversal}
Let $G$ be a digraph, and let $G^{\mathrm{op}}$ be the digraph obtained by reversing every arrow of $G$. With the standard path inner products,
\[
 \Spec\Delta_p^{\loc}(G,T)
 =\Spec\Delta_p^{\loc}(G^{\mathrm{op}},T)
\]
as multisets for any $p\ge0$.
\end{corollary}
\begin{proof}
The signed reversal
\[
 R_p e_{x_0\cdots x_p}
 =(-1)^{p(p+1)/2}e_{x_p\cdots x_0}
\]
is an isometry and satisfies $\partial_pR_p=R_{p-1}\partial_p$. It maps invariant chains in $G$ and $G-T$ onto the corresponding chains in $G^{\mathrm{op}}$ and $G^{\mathrm{op}}-T$, respectively. It therefore induces an isometric chain isomorphism $Q_*(G,T)\to Q_*(G^{\mathrm{op}},T)$. The result follows from Proposition~\ref{prop:isometry}.
\end{proof}

\subsection{Finite-radius determination of the spectrum}\label{sec:locality}

Let $P$ be a path complex. For vertices $x,y\in V(P)$, an undirected path from $x$ to $y$ is a sequence $(x_0,\dots,x_m)$ with $x_0=x$ and $x_m=y$ such that, for each $0\le i<m$, either $x_ix_{i+1}$ or $x_{i+1}x_i$ is an allowed $1$-path. Its length is $m$, and the path $(x)$ has length zero. For a digraph $G$, an undirected path in $P(G)$ is an undirected path in the underlying undirected graph of $G$ connecting $x$ and $y$. The distance $d(x,y)$ between $x$ and $y$ is the length of a shortest undirected path from $x$ to $y$. If no such path exists, we set $d(x,y)=\infty$.

Since the paths are undirected, $d(x,y)=d(y,x)$. One can further verify that $d$ is an extended metric on $V(P)$: it takes values in $[0,\infty]$, vanishes exactly when $x=y$, and satisfies the triangle inequality. The distance $d(x,T)$ between $x$ and a nonempty marked set $T\subseteq V(P)$ is defined by
\[
 d(x,T)=\min_{t\in T}d(x,t).
\]
If there is no undirected path from $x$ to a vertex of $T$, we set $d(x,T)=\infty$. Now, let
\[
 N_s(T)=\{x\in V(P):d(x,T)\le s\},\quad s\ge 0.
\]
The set $N_s(T)$ determines the induced path complex $P|_{N_s(T)}$. For a digraph $G$ and a single marked vertex $v$, write $B_s(v)=G[N_s(\{v\})]$.

\begin{lemma}\label{lem:endpoints}
Let $P$ be a path complex.
For $p\ge1$, write an allowed $p$-chain as $c=\sum_{a,b\in V(P)}c_{a,b}$, where $c_{a,b}$ consists of the terms that start at $a$ and end at $b$. Then $c\in\Omega_p(P)$ if and only if $c_{a,b}\in\Omega_p(P)$ for every pair $a,b$.
\end{lemma}
\begin{proof}
The case $p=1$ follows from $\Omega_1(P)=\A_1(P)$. For $p\ge2$, deleting the first or last vertex of an allowed path gives an allowed face. Thus every non-allowed regular term in $\partial c_{a,b}$ has endpoints $a,b$. Such terms cannot cancel between different endpoint pairs. Hence $\partial c$ is allowed if and only if each $\partial c_{a,b}$ is allowed, proving the assertion.
\end{proof}

\begin{lemma}\label{lem:orthogonal-locality}
For $s\ge p$, we have an isometric isomorphism $Q_p(P|_{N_s(T)},T)\to Q_p(P,T)$.
\end{lemma}
\begin{proof}
Fix $s\ge p$. By Definition~\ref{def:metric}, it suffices to show
\[
 W_p\bigl(P|_{N_s(T)},T\bigr)=W_p(P,T).
\]
For $p=0$, both spaces are spanned by the vertices of $T$. Let $p\ge1$ and $w\in W_p(P,T)$. By Lemma~\ref{lem:endpoints} and the orthogonality of distinct endpoint components, each $w_{a,b}$ also belongs to $W_p(P,T)$. If $d(a,T)+d(b,T)>p$, no allowed $p$-path from $a$ to $b$ meets $T$. Thus $w_{a,b}\in\Omega_p(P-T)$, so $w_{a,b}=0$. For every remaining pair, each vertex $x_i$ of an allowed path $a=x_0\cdots x_p=b$ satisfies
\[
 \begin{aligned}
 d(x_i,T)
 &\le\min\{d(a,T)+i,
              d(b,T)+p-i\}\\
 &\le\frac{d(a,T)+d(b,T)+p}{2}
 \le p.
 \end{aligned}
\]
Hence $w\in\A_p(P|_{N_p(T)})$. Since $w\in\Omega_p(P)$ and $P|_{N_p(T)}$ is induced, $w\in\Omega_p(P|_{N_p(T)})$.

Since $s\ge p$, we have $w\in\Omega_p(P|_{N_s(T)})$. In view of $w\perp\Omega_p(P-T)$, we obtain $w\perp\Omega_p(P|_{N_s(T)}-T)$. It follows that $W_p(P,T)\subseteq W_p(P|_{N_s(T)},T)$. Conversely, write $u\in W_p(P|_{N_s(T)},T)$ as $u=k+w$, where $k\in\Omega_p(P-T)$ and $w\in W_p(P,T)$. Since $u,w\in\Omega_p(P|_{N_s(T)})$, a straightforward calculation shows
\[
 k=u-w\in\Omega_p\bigl(P|_{N_s(T)}\bigr)\cap\Omega_p(P-T)
 =\Omega_p\bigl(P|_{N_s(T)}-T\bigr).
\]
Therefore, one has
\[
 0=\langle u,k\rangle=\|k\|^2,
\]
so $u=w$, proving the equality. The map $[c]\mapsto[c]$ induced by $P|_{N_s(T)}\hookrightarrow P$ identifies the same orthogonal representatives and is an isometric isomorphism.
\end{proof}

\begin{theorem}\label{thm:spectral-locality}
Let $P$ be a path complex and let $T\subseteq V(P)$ be a nonempty marked set. For $p\ge0$ and $s\ge p+1$, the inclusion of pairs $(P|_{N_s(T)},P|_{N_s(T)}-T)\hookrightarrow(P,P-T)$ induces an isometry
\[
 J_p:Q_p\bigl(P|_{N_s(T)},T\bigr)\xrightarrow{\cong}Q_p(P,T).
\]
This map gives the unitary equivalence
\[
 \Delta_p^{\loc}(P,T)
 =J_p\Delta_p^{\loc}\bigl(P|_{N_s(T)},T\bigr)J_p^*.
\]
\end{theorem}
\begin{proof}
By Lemma~\ref{lem:orthogonal-locality}, the inclusion induces isometric isomorphisms in dimensions $p-1$, $p$, and $p+1$ for $p\ge1$, and in dimensions $0$ and $1$ for $p=0$. These maps commute with the local differentials. The argument of Proposition~\ref{prop:isometry} therefore gives
\[
 \Delta_p^{\loc}(P,T)
 =J_p\Delta_p^{\loc}\bigl(P|_{N_s(T)},T\bigr)J_p^*.
\]
Unitary equivalence preserves the spectrum and its multiplicities.
\end{proof}

Taking $s=p+1$ shows that the induced path complex on $N_{p+1}(T)$ determines the local Laplacian in dimension $p$ up to unitary equivalence and hence determines its spectrum.

\begin{corollary}\label{cor:distant-deletion}
Let $G$ be a digraph and let $T\subseteq V(G)$ be a nonempty marked set. For $p\ge0$ and $A\subseteq V(G)\setminus N_{p+1}^G(T)$, the inclusion $G-A\hookrightarrow G$ induces a unitary equivalence
\[
 \Delta_p^{\loc}(G,T)=J_p\Delta_p^{\loc}(G-A,T)J_p^*.
\]
\end{corollary}
\begin{proof}
Every undirected path of length at most $p+1$ starting in $T$ avoids $A$. Thus $G$ and $G-A$ have the same induced path complex $P(G)|_{N_{p+1}(T)}$ on this neighborhood. Theorem~\ref{thm:spectral-locality} gives the unitary equivalences
\[
\Delta_p^{\loc}(G,T)\sim \Delta_p^{\loc}\bigl(P(G)|_{N_{p+1}(T)},T\bigr),\quad \Delta_p^{\loc}(G-A,T)\sim \Delta_p^{\loc}\bigl(P(G)|_{N_{p+1}(T)},T\bigr).
\]
Composing the map from the common neighborhood to $G$ with the inverse of the map to $G-A$ gives $J_p$.
\end{proof}

\begin{remark}
Finite-radius excision for local path homology was established in \cite{LocalHomology}. An isomorphism of chain complexes or of their homology groups alone does not imply unitary equivalence of the local Laplacians.
\end{remark}

\begin{figure}[htbp]
\centering
\includegraphics[width=\textwidth]{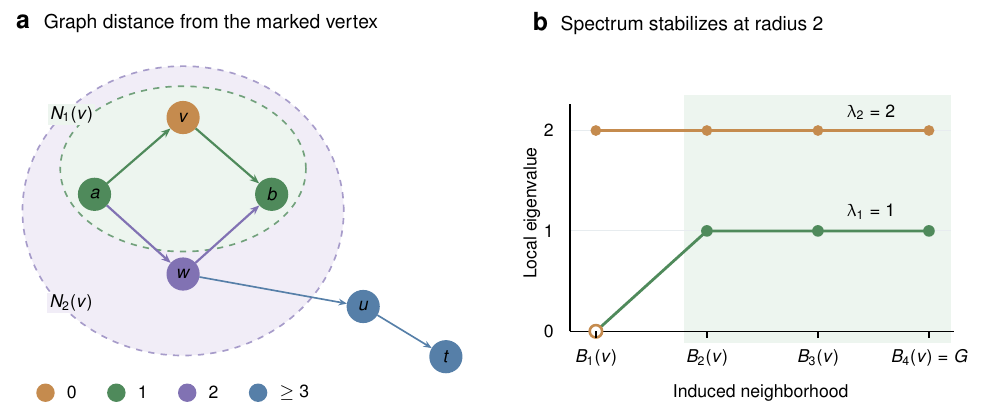}
\caption{\textbf{Finite-radius determination of the full local spectrum.}
\textbf{(a)} A diamond with the additional directed tail $w\to u\to t$. The marked vertex is orange. Green and purple vertices have undirected graph distances one and two from $v$, respectively. Blue vertices have distance at least three. The shaded sets indicate $N_1(v)$ and $N_2(v)$.
\textbf{(b)} The eigenvalues on $B_s(v)$ for $s=1,2,3,4$. The shaded range illustrates Theorem~\ref{thm:spectral-locality} with $p=1$.}
\label{fig:finite-radius}
\end{figure}

Figure~\ref{fig:finite-radius} illustrates the neighborhood bound for the local spectrum in dimension $1$. The first neighborhood contains only the directed interval through $v$. The second includes the alternative path through $w$ and changes the lower eigenvalue from $0$ to $1$. Adding the vertices at distances three and four leaves both eigenvalues unchanged.

For a digraph $G=(V,E)$, let $D=\operatorname{diag}(\deg_G(x))_{x\in V}$ denote the diagonal degree matrix, where $\deg_G(x)=\deg_G^{\mathrm{in}}(x)+\deg_G^{\mathrm{out}}(x)$. The following proposition gives a combinatorial formula for the local Laplacian in dimension $0$.

\begin{proposition}\label{prop:degree-zero}
Let $G=(V,E)$ be a digraph with adjacency matrix $A=(a_{xy})$, and let $T\subseteq V$. In the orthonormal basis $([e_x])_{x\in T}$ of $Q_0(G,T)$, the entries of $L_0^{\loc}(G,T)$ are
\[
 \bigl(L_0^{\loc}(G,T)\bigr)_{xy}
 =\begin{cases}
 \deg_G(x),&x=y,\\
 -(a_{xy}+a_{yx}),&x\ne y,
 \end{cases}
 \quad x,y\in T.
\]
In particular, we have $L_0^{\loc}(G,v)=\bigl(\deg_G(v)\bigr)$.
\end{proposition}
\begin{proof}
Since $\Omega_1(G)=\A_1(G)$, the classes of the edges meeting $T$ form an orthonormal basis of $Q_1(G,T)$. In this basis, the column of $B_1$ for $x\to y$ is the coordinate vector of $[e_y]-[e_x]$, where $[e_z]=0$ for $z\notin T$. Since $d_0=0$, we get
\[
 L_0^{\loc}(G,T)=B_1B_1^{\mathsf T}.
\]
For $x\in T$, each edge incident to $x$ contributes $(\pm1)^2=1$ to $(B_1B_1^{\mathsf T})_{xx}$. An edge $x\to y$ with $x,y\in T$ contributes $(-1)(1)=-1$ to both $(B_1B_1^{\mathsf T})_{xy}$ and $(B_1B_1^{\mathsf T})_{yx}$. Thus the desired formula follows.
\end{proof}

Let $K$ be a finite simplicial complex with ordered vertices. Its associated path complex $P(K)$ consists of the increasing paths $x_0\cdots x_p$ for which $\{x_0,\ldots,x_p\}\in K$ \cite{GLMY}.

For a vertex $v$ of $K$, let $\Delta_p^{\loc}(K,v)$ denote the local simplicial Laplacian on $C_p(K,K-v;\R)$ as defined in \cite{LiuFengLiu}. The following proposition relates this operator to the local path Laplacian.

\begin{proposition}\label{prop:simplicial}
For any $p\ge0$, the local Laplacian $\Delta_p^{\loc}(K,v)$ is unitarily equivalent to the local path Laplacian $\Delta_p^{\loc}(P(K),v)$.
\end{proposition}
\begin{proof}
Each face of an increasing simplex is allowed, so we have $\Omega_p(P(K))=\A_p(P(K))$. Note that the maps
\[
 F_p:C_p(K;\R)\longrightarrow\Omega_p(P(K)),
 \quad [x_0,\ldots,x_p]\longmapsto e_{x_0\cdots x_p}
\]
form an isometric chain isomorphism and send $C_*(K-v;\R)$ onto $\Omega_*(P(K)-v)$. Therefore, we have an isometric chain isomorphism
\[
 \overline F_*:C_*(K,K-v;\R)\longrightarrow Q_*(P(K),v).
\]
By the argument of Proposition~\ref{prop:isometry}, one has
\[
 \Delta_p^{\loc}(P(K),v)
 =\overline F_p\Delta_p^{\loc}(K,v)\overline F_p^*.
\]
This proves the unitary equivalence.
\end{proof}

\section{Examples and computations}\label{sec:examples}

Every example below uses the standard path inner products and the resulting quotient inner products. An edge set consists precisely of the listed edges. Spectra are listed with multiplicity. The spectrum of an operator on the zero vector space is the empty multiset.

\subsection{Basic examples}

\begin{example}\label{ex:isolated}
If $v$ is isolated, then $Q_0(G,v)=\R[e_v]$ and $Q_p(G,v)=0$ for $p\ge1$. Thus, we have
\[
 \Delta_0^{\loc}(G,v)=0,\quad \LH_0(G,v)\cong\R.
\]
There are no local chains in higher dimensions.
\end{example}

\begin{figure}[htbp]
\centering
\includegraphics[width=\textwidth]{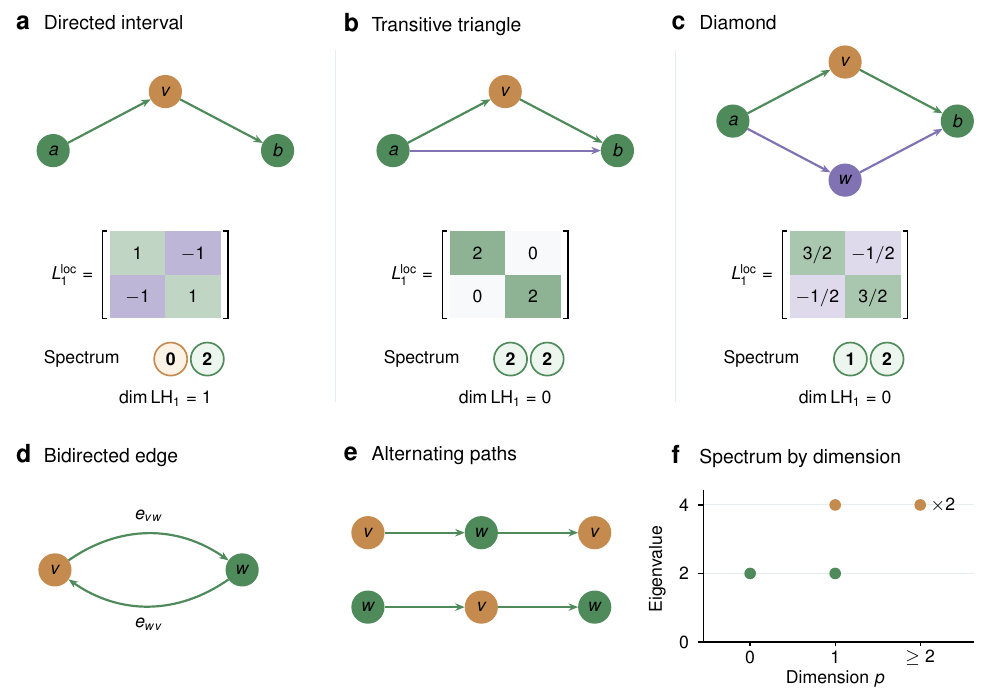}
\caption{\textbf{Local path spectra at a marked vertex in four digraphs.}
\textbf{(a)} The directed interval has no nonzero $\partial$-invariant $2$-chains. Its local spectrum in dimension $1$ is $\{0,2\}$.
\textbf{(b)} Adding $a\to b$ makes $e_{avb}$ $\partial$-invariant and changes this spectrum to $\{2,2\}$.
\textbf{(c)} In the diamond, the non-allowed face $ab$ cancels in $\partial(e_{avb}-e_{awb})$. The class $[e_{avb}-e_{awb}]$ has norm $\sqrt2$, giving the local spectrum $\{1,2\}$.
\textbf{(d,e)} On the bidirected edge, the alternating paths $vwv$ and $wvw$ have the same regular boundary $e_{vw}+e_{wv}$.
\textbf{(f)} The local spectrum is $\{2,4\}$ in dimension $1$ and $\{4,4\}$ in every dimension $p\ge2$.}
\label{fig:motif-spectra}
\end{figure}

\begin{example}\label{ex:interval}
Let $G$ be the interval digraph $a\to v\to b$ in Figure~\ref{fig:motif-spectra}(a). The only allowed $2$-path has boundary
\[
 \partial e_{avb}=e_{vb}-e_{ab}+e_{av}.
\]
Since $ab$ is not allowed, $\Omega_2(G)=0$. So we have
\[
 Q_0=\R[e_v],\quad
 Q_1=\R[e_{av}]\oplus\R[e_{vb}].
\]
In these bases, the local boundary map and the local Laplacians have matrices
\[
 B_1=\begin{pmatrix}1&-1\end{pmatrix},\quad
 L_0^{\loc}=(2),\quad
 L_1^{\loc}=\begin{pmatrix}1&-1\\-1&1\end{pmatrix}.
\]
Therefore, we obtain
\[
 \Spec L_1^{\loc}=\{0,2\},\quad
 \ker L_1^{\loc}=\R\bigl([e_{av}]+[e_{vb}]\bigr),
 \quad \LH_1(G,v)\cong\R.
\]
One can see that the number of zero eigenvalues of $L_1^{\loc}$ is equal to $\dim \LH_1(G,v)=1$.
\end{example}

\begin{example}\label{ex:triangle}
The digraph in Figure~\ref{fig:motif-spectra}(b) is obtained by adding the edge $a\to b$ to the digraph in Example~\ref{ex:interval}. Then $e_{avb}\in\Omega_2(G)$. The deletion $G-v$ contains the edge $a\to b$, so the local spaces in dimensions $0$ and $1$ have the same bases as before, and $Q_2=\R[e_{avb}]$. The local boundary matrices are
\[
 B_1=\begin{pmatrix}1&-1\end{pmatrix},\quad
 B_2=\begin{pmatrix}1\\1\end{pmatrix}.
\]
Consequently, one has
\[
 L_0^{\loc}=(2),\quad
 L_1^{\loc}=\begin{pmatrix}2&0\\0&2\end{pmatrix},\quad
 L_2^{\loc}=(2).
\]
All local homology groups vanish because each local Laplacian has trivial kernel.
\end{example}

\begin{example}\label{ex:diamond}
Let $G$ be the diamond in Figure~\ref{fig:motif-spectra}(c), with vertices $a,v,w,b$ and edges
\[
 a\to v,\quad v\to b,\quad a\to w,\quad w\to b.
\]
The two allowed $2$-paths have the same non-allowed face $ab$. Hence, we have
\[
 \Omega_2(G)=\R z,\quad
 z=e_{avb}-e_{awb},\quad
 \partial z=e_{av}+e_{vb}-e_{aw}-e_{wb}.
\]
The graph $G-v$ is a directed interval, so $\Omega_2(G-v)=0$ and $\|[z]\|^2=2$. The bases for $Q_0$ and $Q_1$ are $[e_v]$ and $([e_{av}],[e_{vb}])$. Using the unit vector $q=[z]/\sqrt2$ in dimension $2$ gives
\[
 B_1=\begin{pmatrix}1&-1\end{pmatrix},\quad
 B_2=\frac1{\sqrt2}\begin{pmatrix}1\\1\end{pmatrix}.
\]
It follows that
\[
 L_0^{\loc}=(2),\quad
 L_1^{\loc}=\begin{pmatrix}\frac32&-\frac12\\-\frac12&\frac32\end{pmatrix},
 \quad L_2^{\loc}=(1),
\]
and $\Spec L_1^{\loc}=\{1,2\}$. All local homology groups vanish.
    
Here $B_1(v)$ is the directed interval $a\to v\to b$, and $B_2(v)=G$. By Example~\ref{ex:interval}, the matrix of the local Laplacian in dimension $1$ at $v$ on $B_1(v)$ is $\left(\begin{smallmatrix}1&-1\\-1&1\end{smallmatrix}\right)$ and has spectrum $\{0,2\}$. On $G$, the spectrum is $\{1,2\}$. Thus the radius $p+1$ in Theorem~\ref{thm:spectral-locality} cannot be reduced to $p$ when $p=1$.
\end{example}

\begin{example}\label{ex:bidirected-edge}
Let $G$ have vertices $v,w$ and arrows $v\to w$ and $w\to v$, with $v$ marked as in Figure~\ref{fig:motif-spectra}(d). Every allowed path alternates between $v$ and $w$. Deleting an interior vertex gives a nonregular face, which vanishes in the regular boundary. Thus the two paths in Figure~\ref{fig:motif-spectra}(e) satisfy
\[
 \partial_2 e_{vwv}=\partial_2 e_{wvw}=e_{vw}+e_{wv}.
\]
The deletion $G-v$ consists of the vertex $w$. Thus $Q_0(G,v)=\R[e_v]$, and in each positive dimension the classes of the two alternating paths form an orthonormal basis of $Q_p(G,v)$. Order these paths by their initial vertex, with $v$ first. Only the first and last faces contribute to the boundary, so
\[
 B_1=\begin{pmatrix}-1&1\end{pmatrix},\qquad
 B_p=\begin{pmatrix}(-1)^p&1\\1&(-1)^p\end{pmatrix}
 \quad (p\ge2).
\]
The local Laplacian matrices are therefore
\[
 L_0^{\loc}=(2),\qquad
 L_1^{\loc}=\begin{pmatrix}3&1\\1&3\end{pmatrix},\qquad
 L_p^{\loc}=4I_2\quad (p\ge2).
\]
In dimension $1$, the chains $[e_{vw}]-[e_{wv}]$ and $[e_{vw}]+[e_{wv}]$ have eigenvalues $2$ and $4$, respectively. All the local eigenvalues in Figure~\ref{fig:motif-spectra}(f) are positive, so $\LH_p(G,v)=0$ for every $p\ge0$.
\end{example}

\subsection{Effects of branching and joint marking}

\begin{figure}[htbp]
\centering
\includegraphics[width=\textwidth]{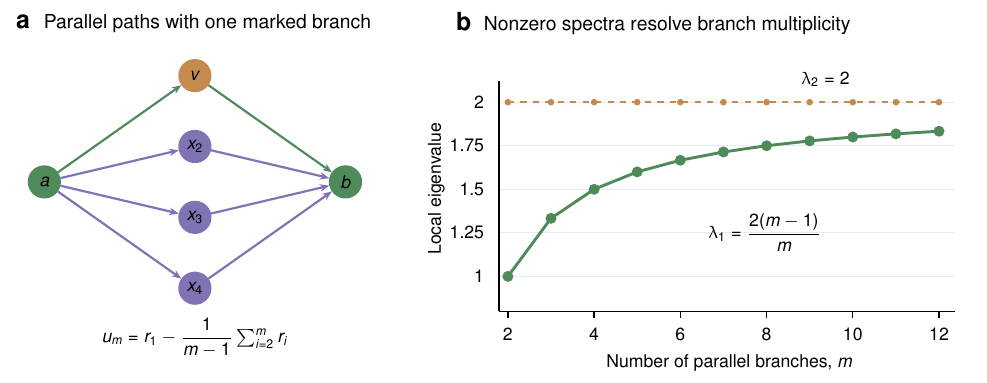}
\caption{\textbf{Parallel paths change nonzero local spectra without changing local homology.}
\textbf{(a)} Four branches of the family $a\to x_i\to b$ are shown, with $v=x_1$ marked. The minimum-norm lift $u_m$ distributes the compensating coefficient equally among the other $m-1$ paths.
\textbf{(b)} The eigenvalues $\lambda_1=2(m-1)/m$ and $\lambda_2=2$ for integer branch counts $m=2,\ldots,12$.}
\label{fig:branch-spectrum}
\end{figure}

\begin{example}\label{ex:branches}
For $m\ge2$, let $G_m$ be the digraph with $m$ branches $a\to x_i\to b$, $1\le i\le m$, as in Figure~\ref{fig:branch-spectrum}.
Mark $v=x_1$ and write $r_i=e_{ax_ib}$. Since $ab$ is not an edge,
\[
 \Omega_2(G_m)=
 \left\{\sum_{i=1}^m c_ir_i:\sum_{i=1}^m c_i=0\right\},
 \quad
 \Omega_2(G_m-v)=
 \left\{\sum_{i=2}^m c_ir_i:\sum_{i=2}^m c_i=0\right\}.
\]
Thus, $Q_2(G_m,v)$ is one dimensional. Its minimum-norm representative with coefficient $1$ on $r_1$ is
\[
 u_m=r_1-\frac1{m-1}\sum_{i=2}^m r_i,
 \quad
 \|u_m\|^2=\frac{m}{m-1}.
\]
Indeed, the chain $u_m$ is $\partial$-invariant and orthogonal to $\Omega_2(G_m-v)$. Set
\[
 q_m=\sqrt{\frac{m-1}{m}}\,[u_m].
\]
With the orthonormal bases $[e_v]$, $([e_{av}],[e_{vb}])$, and $q_m$, the boundary matrices are
\[
 B_1=\begin{pmatrix}1&-1\end{pmatrix},\quad
 B_2=\sqrt{\frac{m-1}{m}}\begin{pmatrix}1\\1\end{pmatrix}.
\]
It follows that
\[
 \begin{gathered}
 L_1^{\loc}(G_m,v)
 =\begin{pmatrix}2-\frac1m&-\frac1m\\-\frac1m&2-\frac1m\end{pmatrix},\\
 \Spec L_1^{\loc}(G_m,v)
 =\left\{\frac{2(m-1)}m,\,2\right\},\quad
 L_2^{\loc}(G_m,v)=\left(\frac{2(m-1)}m\right).
 \end{gathered}
\]
In dimension $0$, we have $L_0^{\loc}=(2)$. These graphs have trivial local homology in every dimension, but their nonzero local spectra vary with $m$. The case $m=2$ is the diamond.
\end{example}

\begin{example}\label{ex:two-marks}
Return to the diamond of Example~\ref{ex:diamond}. First mark the source vertex $a$. Since $\Omega_2(G-a)=0$, use the orthonormal bases $[e_a]$, $([e_{av}],[e_{aw}])$, and $[z]/\sqrt2$ in dimensions $0$, $1$, and $2$, respectively. The boundary matrices are
\[
 B_1^{(a)}=\begin{pmatrix}-1&-1\end{pmatrix},\quad
 B_2^{(a)}=\frac1{\sqrt2}\begin{pmatrix}1\\-1\end{pmatrix}.
\]
It follows that
\[
 L_0^{\loc}(G,a)=(2),\quad
 L_1^{\loc}(G,a)=\frac12\begin{pmatrix}3&1\\1&3\end{pmatrix},\quad
 L_2^{\loc}(G,a)=(1).
\]
Thus $\Spec L_1^{\loc}(G,a)=\{1,2\}$, as for the marked vertex $v$, and all local homology groups vanish.

Now let $T=\{v,w\}$. The deletion consists of the two isolated vertices $a,b$. Use the orthonormal bases
\[
 ([e_v],[e_w]),\quad
 ([e_{av}],[e_{vb}],[e_{aw}],[e_{wb}]),\quad
 [z]/\sqrt2
\]
in dimensions $0$, $1$, and $2$, respectively. A straightforward calculation shows
\[
 B_1=\begin{pmatrix}1&-1&0&0\\0&0&1&-1\end{pmatrix},
 \quad
 B_2=\frac1{\sqrt2}\begin{pmatrix}1\\1\\-1\\-1\end{pmatrix}.
\]
The down Laplacian represented by $B_1^{\mathsf T}B_1$ has eigenvalue $2$ on the chains $[e_{av}]-[e_{vb}]$ and $[e_{aw}]-[e_{wb}]$, while the up Laplacian represented by $B_2B_2^{\mathsf T}$ has eigenvalue $2$ on the chain $[e_{av}]+[e_{vb}]-[e_{aw}]-[e_{wb}]$. Consequently, we obtain
\[
 L_0^{\loc}=2I_2,\quad L_2^{\loc}=(2)
\]
and
\[
 \begin{gathered}
 L_1^{\loc}=B_1^{\mathsf T}B_1+B_2B_2^{\mathsf T}
 =\frac12\begin{pmatrix}
 3&-1&-1&-1\\
 -1&3&-1&-1\\
 -1&-1&3&-1\\
 -1&-1&-1&3
 \end{pmatrix},\quad
 \Spec L_1^{\loc}=\{0,2,2,2\},\\
 \ker L_1^{\loc}
 =\R\bigl([e_{av}]+[e_{vb}]+[e_{aw}]+[e_{wb}]\bigr).
 \end{gathered}
\]
Thus, we have $\LH_1(G,\{v,w\})\cong\R$ and $\LH_1(G,v)=\LH_1(G,w)=0$. The local Laplacian at $\{v,w\}$ is therefore not generally the direct sum of the local Laplacians at $v$ and $w$.
\end{example}

\begin{figure}[htbp]
\centering
\includegraphics[width=\textwidth]{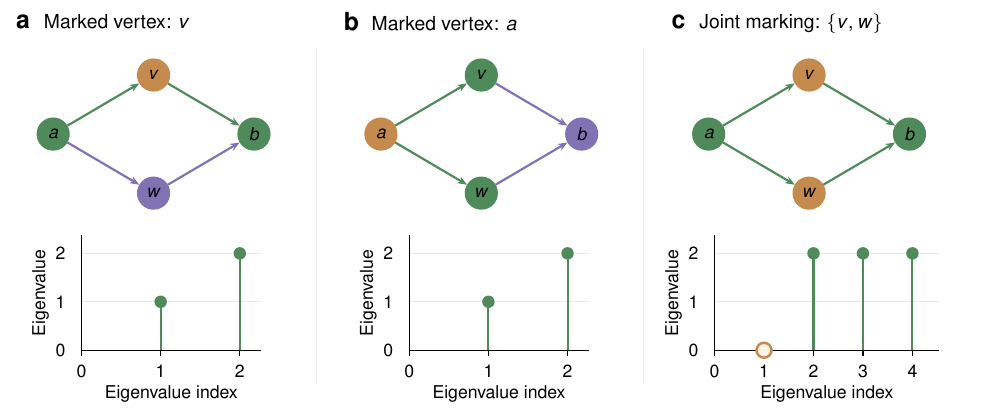}
\caption{\textbf{Local spectra for different marked sets.}
\textbf{(a,b)} Marking $v$ in (a) or $a$ in (b) gives a two-dimensional local $1$-chain space and spectrum $\{1,2\}$.
\textbf{(c)} Marking $\{v,w\}$ gives a four-dimensional local $1$-chain space and spectrum $\{0,2,2,2\}$.}
\label{fig:multiple-marks}
\end{figure}

Figure~\ref{fig:multiple-marks} compares the local spectra for the marked sets $\{v\}$, $\{a\}$, and $\{v,w\}$. Marking $v$ or $a$ gives the same local spectrum $\{1,2\}$ in dimension $1$, although the local boundary matrices differ. Marking both middle vertices produces a four-dimensional local chain space with one harmonic direction. For joint marking, the normalized harmonic chain is $h=\tfrac12([e_{av}]+[e_{vb}]+[e_{aw}]+[e_{wb}])$. This example shows that the local Laplacian at a union of marked sets can differ from the direct sum of the individual operators.

\section{Experiment and application}\label{sec:experiments}

We first examine how well local path features distinguish among isomorphism classes of digraphs. Then we study whether the local spectrum helps identify the same neuron in two animals from its incoming and outgoing synaptic connections. On the same connectomes we measure the time to compute the global spectra and the local spectra at all neurons.

\subsection{Local spectral discrimination}\label{sec:five-vertex-spectra}

Let $V=\{a,b,c,d,e\}$. For each unordered pair of vertices, we either omit the edge or choose one of its two directions. This gives $3^{\binom{5}{2}}=59{,}049$ labeled digraphs, including the empty graph. By Burnside's lemma, these digraphs fall into exactly $582$ isomorphism classes under vertex relabeling.

For $v\in V$, let $S_p(G,v)=\Spec\Delta_p^{\loc}(G,v)$ be the local spectrum in dimension $p$. In this section, we consider the cases $p=0$ and $p=1$. To compare digraphs up to isomorphism, we use the multiset
\[
 \mathcal S(G)=\bigl\{\!\bigl\{\,(S_0(G,v),S_1(G,v))\mid v\in V\,\bigr\}\!\bigr\}.
\]
An isomorphism maps each marked vertex to the corresponding marked vertex and preserves both spectra. Let $\omega_G(v)$ denote the total degree of $v$, the sum of its in-degree and out-degree. By Proposition~\ref{prop:degree-zero}, $S_0(G,v)=\{\omega_G(v)\}$. Since the incident arrows form a basis of $Q_1(G,v)$, the cardinality of $S_1(G,v)$ is also $\omega_G(v)$. Thus $S_1(G,v)$ determines $S_0(G,v)$.

For each of the 582 representatives, we computed the local Laplacian matrices in dimensions $0$ and $1$ at all five vertices and verified them against the quotient construction in Section~\ref{sec:matrices}. Spectra were compared via their characteristic polynomials. Relabeling yields the spectra for all $295{,}245$ labeled graph–vertex pairs. With labels fixed, recording the local spectra at the first $k$ vertices in the order $a,b,c,d,e$ produces $78$, $657$, $2390$, $5225$, and $5{,}449$ distinct ordered lists for $k=1,\dots,5$, respectively.

Table~\ref{tab:five-vertex-spectra} compares four choices of vertex data. Write $\beta_p=\beta_p(G,v)$ and $\omega=\omega_G(v)$. Each choice partitions the $582$ isomorphism types into classes according to the multiset of the five vertex measurements. If the class sizes are $n_1,\ldots,n_r$, then $\sum_{i=1}^r\binom{n_i}{2}$ of the $\binom{582}{2}=169{,}071$ unordered pairs have identical data.

\begin{table}[htbp]
\centering
\caption{Distinguishing five-vertex digraphs without reciprocal edges. Each isomorphism type is counted once. A singleton is a class containing one type.}
\label{tab:five-vertex-spectra}
\small
\begin{tabular*}{0.8\linewidth}{@{\extracolsep{\fill}}lrrrr@{}}
\hline\noalign{\vskip 3pt}
Data at $v$ & Classes & Singletons & Pairs distinguished & Proportion\\
\noalign{\vskip 2pt}\hline\noalign{\vskip 2pt}
$(\beta_0,\beta_1)$      &  28 &  5 & $146{,}022$ & $86.37\%$\\
\quad $\omega$                   &  31 &  3 & $158{,}370$ & $93.67\%$\\
$(\omega,\beta_0,\beta_1)$ &  80 & 14 & $165{,}365$ & $97.81\%$\\
$(S_0,S_1)$              & 142 & 40 & $167{,}082$ & $98.82\%$\\
\noalign{\vskip 2pt}\hline
\end{tabular*}
\end{table}

\begin{figure}[htbp]
\centering
\includegraphics[width=\textwidth]{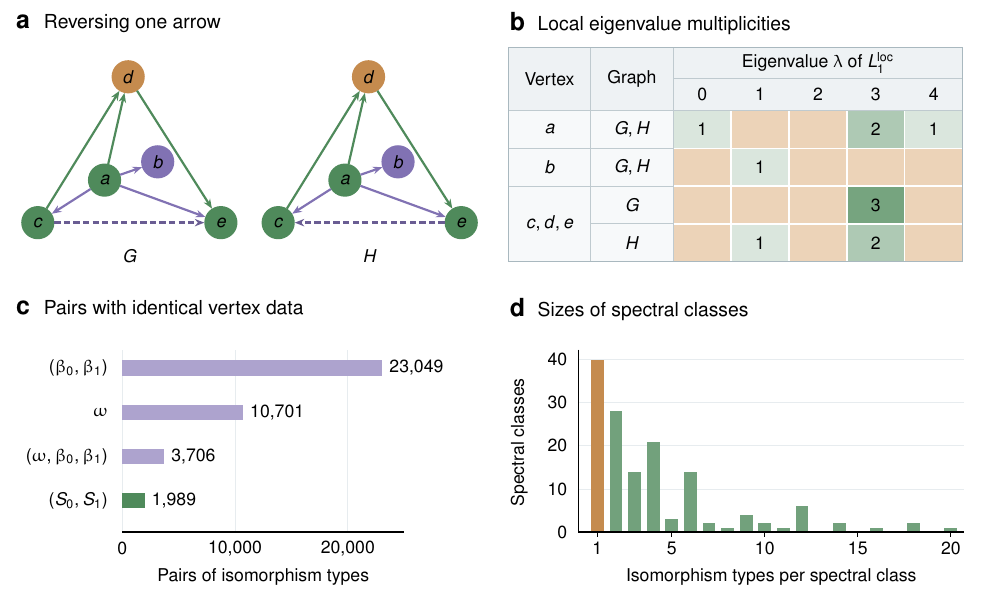}
\caption{\textbf{Local spectra of five-vertex digraphs without reciprocal edges.}
\textbf{(a)} The dashed arrow is reversed between $G$ and $H$. The marked vertex $d$ is orange.
\textbf{(b)} Each column gives the multiplicity of $\lambda$ in the local spectrum of $L_1^{\loc}$ at the indicated vertices. A blank pale orange cell means multiplicity zero. Green cells give positive multiplicities, with darker shades for larger values. For each $v\in\{c,d,e\}$, ordering the eigenvalues as $\lambda_1\leq\lambda_2\leq\lambda_3$ gives $\lambda_1(G,v)=3$ and $\lambda_1(H,v)=1$.
\textbf{(c)} Counts of indistinguishable pairs for the four choices of vertex data in Table~\ref{tab:five-vertex-spectra}. Here $\omega$ is the total vertex degree.
\textbf{(d)} Numbers of spectral classes of each size. The bar at $k$ counts the spectral classes containing exactly $k$ isomorphism types. For example, the orange bar shows that 40 types are uniquely identified by their local spectra.}
\label{fig:five-vertex-triangle}
\end{figure}

Among the $10{,}701$ pairs of nonisomorphic graphs with the same multiset of total vertex degrees, the local spectra in dimension $1$ distinguish $8{,}712$ pairs ($81.41\%$). In $3{,}706$ pairs, both graphs also have the same multiset of triples $(\omega,\beta_0,\beta_1)$. The positive eigenvalues distinguish $1{,}717$ of these pairs ($46.33\%$). Over the entire enumeration, $\mathcal S$ distinguishes $167{,}082$ of the $169{,}071$ pairs ($98.82\%$). Of the $142$ spectral classes, $40$ are singletons, so the spectra uniquely identify $40$ of the $582$ types ($6.87\%$). The largest spectral class contains $20$ types. Figure~\ref{fig:five-vertex-triangle}(c,d) shows the pair counts and the distribution of class sizes.

Figure~\ref{fig:five-vertex-triangle}(a,b) shows one of the $1{,}717$ pairs distinguished by the positive eigenvalues. The graph $G$ has the four arrows from $a$ to $b,c,d,e$ and the three arrows $c\to d$, $d\to e$, and $c\to e$. It is acyclic. Reversing $c\to e$ gives $H$ and creates the directed cycle $c\to d\to e\to c$. This reversal preserves the total degree and both local Betti numbers at each vertex.

Mark $d$ and order the incident edges as $ad,cd,de$. In both graphs $B_1=(1,1,-1)$, and the local boundaries of $e_{acd}$ and $e_{ade}$ have coordinates $(-1,1,0)^{\mathsf T}$ and $(1,0,1)^{\mathsf T}$, respectively. In $G$ the chain $e_{cde}$ contributes a third column $u=(0,1,1)^{\mathsf T}$ to $B_2$. In $H$, each $2$-path around the directed triangle has a non-allowed face that occurs in the boundary of no other allowed $2$-path. Thus only the two common columns contribute, and
\[
 L_1^{\loc}(G,d)=3I_3,\qquad
 L_1^{\loc}(H,d)=3I_3-uu^{\mathsf T}
 =\begin{pmatrix}3&0&0\\0&2&-1\\0&-1&2\end{pmatrix}.
\]
Since $\|u\|^2=2$, reversing this arrow changes one eigenvalue from $3$ to $1$ at $d$. The same change occurs at $c$ and $e$. The common spectra in dimension $1$ at $a$ and $b$ are $\{0,3,3,4\}$ and $\{1\}$, respectively. In the $\lambda=0$ column of Figure~\ref{fig:five-vertex-triangle}(b), the green entry $1$ at $a$ gives $\beta_1(G,a)=\beta_1(H,a)=1$. The pale orange cell at $b$ gives $\beta_1(G,b)=\beta_1(H,b)=0$.

By Corollary~\ref{cor:reversal}, reversing all arrows preserves the local spectra at every vertex. Of the $1{,}989$ pairs with the same $\mathcal S$, exactly $240$ are related by reversal of all arrows.

\subsection{Matching neurons across adult connectomes}\label{sec:adult-neuron-matching}

We used the WormWiring connection tables for specimens~7 and~8, the two adult \emph{C.~elegans} connectomes reconstructed by Witvliet et al.~\cite{WitvlietConnectomes}. In each table we retained the named sensory neurons, interneurons, motor neurons and modulatory neurons. This gave the same $180$ neuron names in both animals. We drew an arrow $u\to v$ whenever the table recorded at least one chemical synapse from $u$ to a distinct neuron $v$. Each graph has $1{,}933$ arrows, with reciprocal connections retained as two arrows.

For each neuron in specimen~7, we ranked all $180$ neurons in specimen~8 by increasing distance and recorded the rank of the neuron with the same name. Repeating this procedure in the opposite direction gave $360$ rankings, two for each of the $180$ neuron identities. The published names define the cohort and identify the correct match. The ranking uses only the graph structure.

Write $\mathbf d(u)=(d_{\mathrm{in}}(u),d_{\mathrm{out}}(u))$ for the degree pair at $u$. By Theorem~\ref{thm:spectral-locality}, the spectrum $\Spec L_1^{\loc}(G,u)$ can be computed on $B_2(u)$, the subgraph induced by the vertices reachable from $u$ in at most two steps of the underlying undirected graph. We verified this equality at all $360$ marked neurons by computing the spectra on the full graphs. Restricting to $B_1(u)$ changed at least one eigenvalue at every mark by more than $10^{-9}$. Let $F_u(t)$ be the fraction of eigenvalues in $\Spec L_1^{\loc}(G,u)$ that are at most $t$, counted with multiplicity. For neurons $u$ and $v$ in different animals, we use the degree distance and the $1$-Wasserstein distance between their spectral distributions,
\[
 D_{\mathrm{deg}}(u,v)=\|\mathbf d(u)-\mathbf d(v)\|_1,\qquad
 D_{\mathrm{spec}}(u,v)=\int_0^\infty |F_u(t)-F_v(t)|\,\mathrm{d}t.
\]
For the combined ranking we use $D_{\mathrm{deg}}/m_{\mathrm{deg}}+D_{\mathrm{spec}}/m_{\mathrm{spec}}$, where $m_{\mathrm{deg}}=14$ and $m_{\mathrm{spec}}\approx5.266$ are the medians of the corresponding distances over all $180^2$ cross-animal candidate pairs.

\begin{figure}[p]
\centering
\includegraphics[width=\textwidth]{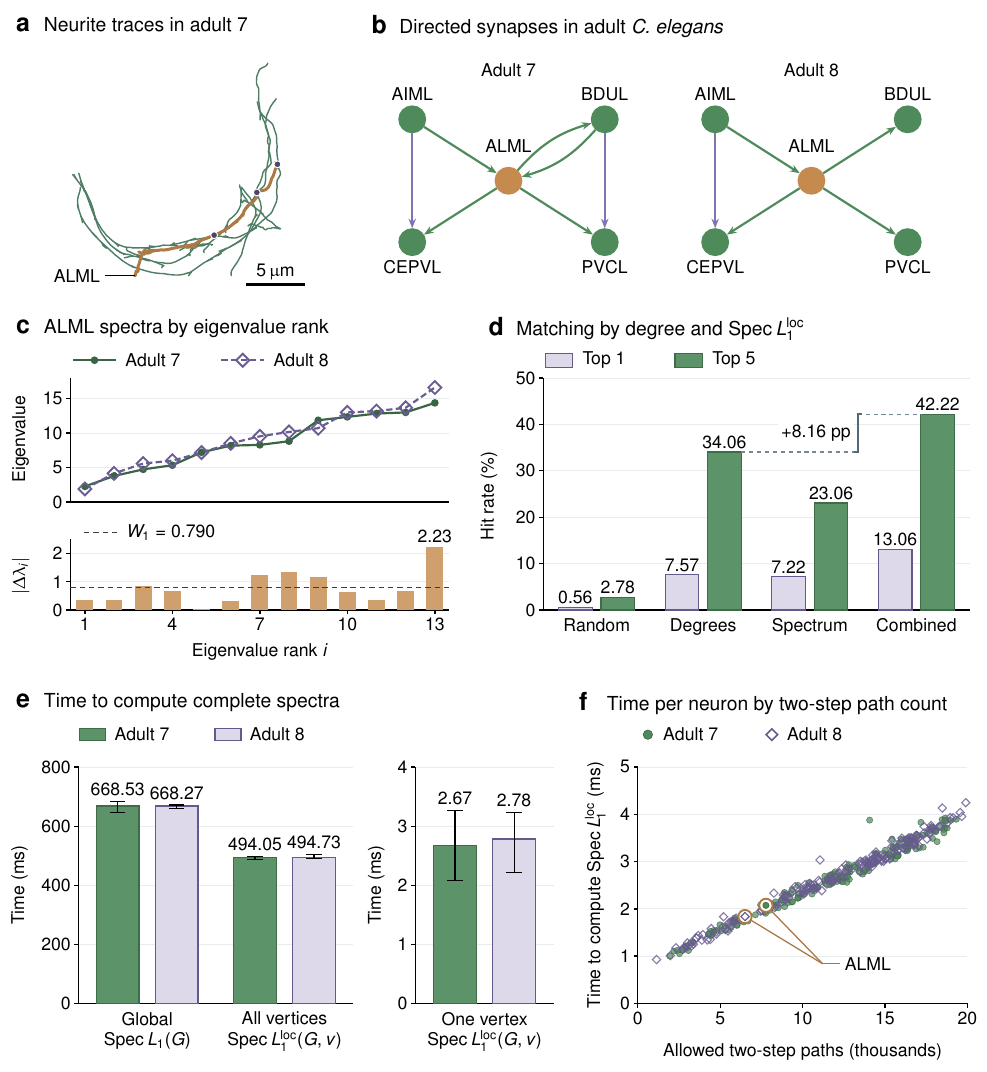}
\caption{\textbf{Neuron matching and spectral computation in adult \emph{C.~elegans} connectomes.}
\textbf{(a)} Projected neurite traces from specimen~7 \cite{WitvlietConnectomes,WitvlietTraceData}. Purple dots on ALML (orange) mark one input from AIML and two from BDUL.
\textbf{(b)} Chemical synapses among these five neurons in both adults. Arrows point from the presynaptic to the postsynaptic neuron.
\textbf{(c)} ALML eigenvalues paired by rank. Bars show the 13 absolute differences $|\lambda_i^{(7)}-\lambda_i^{(8)}|$. Their mean (dashed line) is the $1$-Wasserstein distance.
\textbf{(d)} Top-one and top-five hit rates in both directions. The bracket marks the top-five gain from adding local spectra to degrees.
\textbf{(e)} Bars and whiskers give medians and interquartile ranges. Global and all-vertex times summarize ten runs. The single-vertex bars summarize the $180$ per-neuron medians on a separate scale.
\textbf{(f)} Each point gives the ten-round median time to compute $\Spec L_1^{\loc}(G,v)$ for one neuron $v$ versus the number of allowed two-step paths in $B_2(v)$. Orange rings mark ALML.}
\label{fig:adult-neuron-matching}
\end{figure}

The top-$k$ hit rate is the fraction of rankings in which the correct name appears among the first $k$ candidates, averaged over uniform orderings within each tie. A uniformly random ranking gives an expected hit rate of $k/180$. Table~\ref{tab:adult-neuron-matching} and Figure~\ref{fig:adult-neuron-matching}(d) give the results. For this pair of animals, adding $\Spec L_1^{\loc}$ to the degrees raised the top-five hit rate by $8.16$ percentage points. Of the $360$ rankings, the combined method placed the correct name first in $47$ and among the first five in $152$.

\begin{table}[htbp]
\centering
\caption{Exact-name neuron matching between adult specimens~7 and~8. Each ranking includes all $180$ candidate neurons in the other animal. The highest hit rate in each column is in bold.}
\label{tab:adult-neuron-matching}
\small
\begin{tabular*}{0.85\linewidth}{@{\extracolsep{\fill}}lrr@{}}
\hline\noalign{\vskip 3pt}
Ranking method & Top-1 hit rate & Top-5 hit rate\\
\noalign{\vskip 2pt}\hline\noalign{\vskip 2pt}
Random ranking                     & $0.56\%$  & $2.78\%$\\
In- and out-degrees                & $7.57\%$  & $34.06\%$\\
$\Spec L_1^{\loc}$                 & $7.22\%$  & $23.06\%$\\
Degrees and $\Spec L_1^{\loc}$     & \textbf{13.06\%} & \textbf{42.22\%}\\
\noalign{\vskip 2pt}\hline
\end{tabular*}
\end{table}

The neuron ALML illustrates how the spectrum resolves a degree tie. ALML in specimen~7 has in-degree $4$ and out-degree $9$. In specimen~8, ALML, SDQL and ASIR all have this degree pair. Their spectral distances from ALML in specimen~7 are $0.790$, $1.052$ and $3.946$, respectively. The combined score ranks ALML first among all $180$ candidates. Figure~\ref{fig:adult-neuron-matching}(a,b) locates three synaptic inputs along ALML in specimen~7 and shows the connections among the same five neurons in each animal. Panel~(c) compares the two local spectra at ALML in ascending eigenvalue order.

\subsection{Cost of local and global path spectra}\label{sec:adult-laplacian-timing}

\paragraph{Complexity in dimension~$1$.}
After removing isolated vertices, let $n$ and $N$ denote the numbers of vertices and edges of the simple loopless digraph $G$. Set $d_v=d_{\mathrm{in}}(v)+d_{\mathrm{out}}(v)$ and $D=\max_v d_v$. The global edge basis has $N$ elements. At $v$, the local basis consists of the $d_v$ incident edges. Therefore $L_1(G)$ and $L_1^{\loc}(G,v)$ have sizes $N$ and $d_v$, respectively.

Let $k_{ac}$ be the number of two-step paths from $a$ to $c$, and let $R=\sum_{a,c}k_{ac}$. Enumerating these paths from the dense adjacency matrix takes $O(n^2+R)$ time. The down Laplacian takes $O(\sum_v d_v^2)=O(N^2)$ time to assemble, since $\sum_v d_v=2N$. For the up Laplacian, each boundary column has at most three nonzero entries, forming their outer products takes $O(R)$ time. Sorting the paths by endpoint pair takes $O(R\log(2+R))$ time. If $a\ne c$, $k_{ac}>0$, and the edge $a\to c$ is absent, then in any $\partial$-invariant $2$-chain the coefficients of the two-step paths from $a$ to $c$ sum to zero. The projector $I-\mathbf1\mathbf1^{\mathsf T}/k_{ac}$ therefore gives a rank-one correction on at most $2k_{ac}$ edges, at cost $O(k_{ac}^2)$. Since the intermediate vertices are distinct,
\[
 k_{ac}\le\min\{d_{\mathrm{out}}(a),d_{\mathrm{in}}(c)\},
 \quad
 R\le\sum_{a,c}k_{ac}^2
 \le\sum_{a,c}d_{\mathrm{out}}(a)d_{\mathrm{in}}(c)=N^2.
\]
Thus, using $n\le 2N$, global assembly takes $O(N^2+R\log(2+R))$ time and $O(N^2)$ space. With dense diagonalization, the total cost is $O(N^3)$ time and $O(N^2)$ space.

Let $r$ bound the total number of paths of lengths zero, one, and two in each $B_2(v)$. Scanning the adjacency rows and columns of $v$ and its neighbors identifies the vertices of $B_2(v)$ in $O(n(1+d_v))$ time and space. The assembly estimate above gives $O(r^2)$ time and space for the local matrix, including the induced adjacency matrix. Since the local edge basis has $d_v\le r$ elements, dense diagonalization gives a total cost of $O(n(1+d_v)+r^2+d_v^3)$ time and $O(n(1+d_v)+r^2)$ working space.

\begin{table}[htbp]
\centering
\caption{Cost of computing complete path spectra in dimension~$1$. Working storage is in addition to the input adjacency matrix.}
\label{tab:dimension-one-cost}
\small
\begin{tabular*}{0.97\linewidth}{@{\extracolsep{\fill}}lcc@{}}
\hline
Computation & Time & Peak working storage\\
\hline
Global spectrum & $O(N^3)$ & $O(N^2)$\\
Local spectrum at $v$ & $O(n(1+d_v)+r^2+d_v^3)$ & $O(n(1+d_v)+r^2)$\\
Local spectra at all vertices & $O(nN+Nr^2)$ & $O(n(1+D)+r^2)$\\
\hline
\end{tabular*}
\end{table}

For the sequential pass in Table~\ref{tab:dimension-one-cost}, sum the three single-vertex time bounds using $n\le 2N$ and $D\le r$:
\[
 \sum_v n(1+d_v)=n(n+2N),\qquad
 \sum_v r^2=nr^2,\qquad
 \sum_v d_v^3\le D^2\sum_v d_v=2ND^2.
\]
The total time is therefore $O(nN+Nr^2)$. The pass reuses the single-vertex workspace and stores $\sum_v d_v=2N$ eigenvalues. Since $2N\le nD$, peak storage is $O(n(1+D)+r^2)$. If $r=o(N)$, then $(nN+Nr^2)/N^3=n/N^2+(r/N)^2\to0$. Computing all local spectra is therefore asymptotically faster than dense global diagonalization.

\paragraph{Adult connectome benchmark.}
Each adult connectome has $180$ vertices and $1{,}933$ directed edges, with $R=20{,}169$ and $20{,}504$ allowed two-step paths in specimens~7 and~8, respectively. We compared the time to compute the global spectrum $\Spec L_1(G)$, the local spectrum $\Spec L_1^{\loc}(G,v)$ at one vertex, and the full family $\bigl(\Spec L_1^{\loc}(G,v)\bigr)_{v\in V(G)}$. By Theorem~\ref{thm:spectral-locality}, each local spectrum can be computed independently on $B_2(v)$.

Timing started with the adjacency matrix in memory and ended after the last eigenvalue was computed. For each local spectrum, we included the time required to extract its two-step neighborhood. The computations were performed on a Windows~11 laptop with an Intel Core Ultra~9 processor and $32$~GB of RAM. All matrix operations used double precision on a single CPU thread. After one warm-up run, we recorded ten rounds, timing one global computation and one complete sequential pass over all vertices in each round. We also recorded the time for each vertex in the pass. For the single-vertex result, we first took the median over the ten rounds for each vertex and then the median over all $180$ vertices.

\begin{table}[htbp]
\centering
\caption{Median elapsed time in milliseconds for path spectra in dimension~$1$. The all-vertices entry gives the measured time for one complete pass over the $180$ vertices.}
\label{tab:adult-laplacian-timing}
\small
\begin{tabular*}{0.82\linewidth}{@{\extracolsep{\fill}}lrr@{}}
\hline\noalign{\vskip 2pt}
Computation & Adult 7 & Adult 8\\
\noalign{\vskip 2pt}\hline\noalign{\vskip 2pt}
Global spectrum, $\Spec L_1(G)$ & 668.53 & 668.27\\
One vertex, $\Spec L_1^{\loc}(G,v)$ & 2.67 & 2.78\\
All vertices, $\bigl(\Spec L_1^{\loc}(G,v)\bigr)_{v\in V(G)}$ & 494.05 & 494.73\\
\noalign{\vskip 2pt}\hline
\end{tabular*}
\end{table}

\begin{samepage}
The median time for one local spectrum was under $3$~ms in both animals (Table~\ref{tab:adult-laplacian-timing} and Figure~\ref{fig:adult-neuron-matching}(e)). In each animal the complete pass over all $180$ vertices took about $26\%$ less time than computing the global spectrum.
\par
\end{samepage}

Figure~\ref{fig:adult-neuron-matching}(f) shows how the time to compute a neuron's local spectrum varies with the number $r_v$ of allowed two-step paths in $B_2(v)$. The Spearman rank correlations were $0.984$ in specimen~7 and $0.990$ in specimen~8. For ALML, $d_v=13$ in both animals. The two neighborhoods contained $7{,}774$ and $6{,}506$ allowed two-step paths, respectively. The corresponding median times to compute the local spectrum were $2.07$ and $1.84$~ms.


\begin{thebibliography}{99}

\bibitem{BattistonHigherOrder}
F.~Battiston, G.~Cencetti, I.~Iacopini, V.~Latora, M.~Lucas, A.~Patania, J.-G.~Young and G.~Petri.
\newblock Networks beyond pairwise interactions: structure and dynamics.
\newblock \emph{Physics Reports}, 874:1--92, 2020.
\newblock \url{https://doi.org/10.1016/j.physrep.2020.05.004}.

\bibitem{GodsilRoyle}
C.~Godsil and G.~Royle.
\newblock \emph{Algebraic Graph Theory}.
\newblock Graduate Texts in Mathematics, vol.~207. Springer, New York, 2001.
\newblock \url{https://doi.org/10.1007/978-1-4613-0163-9}.

\bibitem{ReimannCliques}
M.~W.~Reimann, M.~Nolte, M.~Scolamiero, K.~Turner, R.~Perin, G.~Chindemi, P.~D{\l}otko, R.~Levi, K.~Hess and H.~Markram.
\newblock Cliques of neurons bound into cavities provide a missing link between structure and function.
\newblock \emph{Frontiers in Computational Neuroscience}, 11:48, 2017.
\newblock \url{https://doi.org/10.3389/fncom.2017.00048}.

\bibitem{JiangHodgeRanking}
X.~Jiang, L.-H.~Lim, Y.~Yao and Y.~Ye.
\newblock Statistical ranking and combinatorial Hodge theory.
\newblock \emph{Mathematical Programming}, 127:203--244, 2011.
\newblock \url{https://doi.org/10.1007/s10107-010-0419-x}.

\bibitem{GLMY}
A.~Grigor'yan, Y.~Lin, Y.~Muranov and S.-T.~Yau.
\newblock \emph{Homologies of path complexes and digraphs}.
\newblock arXiv:1207.2834, 2012.
\newblock \url{https://arxiv.org/abs/1207.2834}.

\bibitem{GJMY}
A.~Grigor'yan, R.~Jimenez, Y.~Muranov and S.-T.~Yau.
\newblock On the path homology theory of digraphs and Eilenberg--Steenrod axioms.
\newblock \emph{Homology, Homotopy and Applications}, 20(2):179--205, 2018.

\bibitem{ChenLiuPathTopology}
D.~Chen, J.~Liu, J.~Wu, G.-W.~Wei, F.~Pan and S.-T.~Yau.
\newblock Path topology in molecular and materials sciences.
\newblock \emph{The Journal of Physical Chemistry Letters}, 14(4):954--964, 2023.
\newblock \url{https://doi.org/10.1021/acs.jpclett.2c03706}.

\bibitem{WangWeiPathLaplacian}
R.~Wang and G.-W.~Wei.
\newblock Persistent path Laplacian.
\newblock \emph{Foundations of Data Science}, 5(1):26--55, 2023.
\newblock \url{https://doi.org/10.3934/fods.2022015}.

\bibitem{ChowdhuryMemoliPath}
S.~Chowdhury and F.~M\'emoli.
\newblock Persistent path homology of directed networks.
\newblock In \emph{Proceedings of the Twenty-Ninth Annual ACM--SIAM Symposium on Discrete Algorithms}, pages 1152--1169, 2018.
\newblock \url{https://doi.org/10.1137/1.9781611975031.75}.

\bibitem{DeyLiWangPath}
T.~K.~Dey, T.~Li and Y.~Wang.
\newblock An efficient algorithm for 1-dimensional (persistent) path homology.
\newblock \emph{Discrete \& Computational Geometry}, 68:1102--1132, 2022.
\newblock \url{https://doi.org/10.1007/s00454-022-00430-8}.

\bibitem{KazemiDynamicGraphs}
S.~M.~Kazemi, R.~Goel, K.~Jain, I.~Kobyzev, A.~Sethi, P.~Forsyth and P.~Poupart.
\newblock Representation learning for dynamic graphs: a survey.
\newblock \emph{Journal of Machine Learning Research}, 21(70):1--73, 2020.
\newblock \url{https://www.jmlr.org/papers/v21/19-447.html}.

\bibitem{LocalHomology}
J.~Liu and K.~Liu.
\newblock \emph{On the local homology of path complexes}.
\newblock Manuscript, 2026.

\bibitem{LiuFengLiu}
J.~Liu, H.~Feng and K.~Liu.
\newblock \emph{Local Laplacian: theory and models for data analysis}.
\newblock arXiv:2603.07591, 2026.
\newblock \url{https://arxiv.org/abs/2603.07591}.

\bibitem{ChenWangLin2026}
Y.~Chen, S.~Wang and H.~Lin.
\newblock Topology structure optimization of reservoirs using GLMY homology.
\newblock \emph{Neural Networks}, 109004, 2026.

\bibitem{WitvlietConnectomes}
D.~Witvliet, B.~Mulcahy, J.~K.~Mitchell, Y.~Meirovitch, D.~R.~Berger, Y.~Wu, Y.~Liu, W.~X.~Koh, R.~Parvathala, D.~Holmyard, R.~L.~Schalek, N.~Shavit, A.~D.~Chisholm, J.~W.~Lichtman, A.~D.~T.~Samuel and M.~Zhen.
\newblock Connectomes across development reveal principles of brain maturation.
\newblock \emph{Nature}, 596:257--261, 2021.
\newblock \url{https://doi.org/10.1038/s41586-021-03778-8}.

\bibitem{WitvlietTraceData}
D.~Witvliet et al.
\newblock Neurite traces and chemical synapse annotations for specimen~7.
\newblock Data accompanying \emph{Connectomes across development reveal principles of brain maturation}, 2021.
\newblock \href{https://github.com/dwitvliet/nature2021/tree/0646af9d25896ae660f97d462eab2d67282f5625/data}{Archived repository version}.

\end{thebibliography}
\end{document}